\documentclass[11pt,leqno]{amsart}
\usepackage{amssymb,verbatim,enumerate,ifthen}
\usepackage[mathscr]{eucal}
\usepackage[utf8]{inputenc}
\usepackage[T1]{fontenc}
\usepackage{cite}
\def\H{\mathbb{H}}
\def\N{\mathbb{N}}
\def\R{\mathbb{R}}

\def\P{\mathscr{P}}

\def\L{\mathscr{L}}
\def\G{\mathscr{G}}
\def\A{\mathscr{A}}
\def\H{\mathscr{H}}

\newtheorem{theorem}{Theorem}[section]
\newtheorem*{theorem*}{Theorem}
\def\Thm#1#2{\ifthenelse{\equal{#1}{*}}{\begin{theorem*}#2\end{theorem*}}
             {\begin{theorem}\label{T#1}#2\end{theorem}}}
\newtheorem{Atheorem}{Theorem}

\def\thm#1{Theorem~\ref{T#1}}
\newtheorem{proposition}[theorem]{Proposition}
\newtheorem*{proposition*}{Proposition}
\def\Prp#1#2{\ifthenelse{\equal{#1}{*}}{\begin{proposition*}#2\end{proposition*}}
{\begin{proposition}\label{P#1}#2\end{proposition}}}

\newtheorem{corollary}[theorem]{Corollary}
\newtheorem*{corollary*}{Corollary}
\def\Cor#1#2{\ifthenelse{\equal{#1}{*}}{\begin{corollary*}#2\end{corollary*}}
             {\begin{corollary}\label{C#1}#2\end{corollary}}}

\newtheorem{lemma}[theorem]{Lemma}
\newtheorem*{lemma*}{Lemma}
\def\Lem#1#2{\ifthenelse{\equal{#1}{*}}{\begin{lemma*}#2\end{lemma*}}
             {\begin{lemma}\label{L#1}#2\end{lemma}}}

\theoremstyle{definition}
\newtheorem{remark}[theorem]{Remark}
\newtheorem*{remark*}{Remark}
\def\Rem#1#2{\ifthenelse{\equal{#1}{*}}{\begin{remark}\rm #2\end{remark}}
             {\begin{remark}\label{R#1}\rm #2\end{remark}}}

\newtheorem{example}[theorem]{Example}
\newtheorem*{example*}{Example}
\def\Exa#1#2{\ifthenelse{\equal{#1}{*}}{\begin{example*}\rm #2\end{example*}}
             {\begin{example}\label{Ex#1}\rm #2\end{example}}}

\def\eq#1{{\rm(\ref{E#1})}}
\def\Eq#1#2{\ifthenelse{\equal{#1}{*}}
  {\begin{equation*}\begin{aligned}#2\end{aligned}\end{equation*}}
  {\begin{equation}\begin{aligned}\label{E#1}#2\end{aligned}\end{equation}}}

\begin{document}
\begin{flushright}
\end{flushright}
\vspace{5mm}

\date{\today}

\title{Revisiting Ostrowski's Inequality}

\author[A. R. Goswami]{Angshuman R. Goswami}
\address[A. R. Goswami]{Department of Mathematics, University of Pannonia, 
H-8200 Veszprem, Hungary}
\email{goswami.angshuman.robin@mik.uni-pannon.hu}
\thanks{The author’s research was supported by the EKÖP Scholarship (2024-2.1.1-EKÖP-2024-00025/29)\,, funded by the Research Fellowship Programme of the Ministry of Culture and Innovation, Government of Hungary.}

\begin{abstract} The main objective of this paper is to present Ostrowski's inequality for a broader class of functions and to propose a refinement to the classical version of it. The original Ostrowski's inequality can be stated as follows\\

"If $f:[a,b]\to\R$ is differentiable and $f'\in L^{\infty}[a, b]$, then for any $p\in\,]a,b[\,$, the following functional inequality holds:
\Eq{*}{
\Bigg|f(p)-\dfrac{1}{b-a}\int_{a}^{b}f(t)\,dt\Bigg|\leq \dfrac{(p-a)^2+(b-p)^2}{2(b-a)}\Big\| f'\Big\|_a^b\,.\,^{^{^{^{^"}}}}
}
We relax the condition of differentiability and show that even if $f\in C[a,b]$ is non-differentiable at the points $p_{_{1}},\cdots,p_{_{n}}$, then for any $p\in\,]a,b[\,\setminus\overset{n}{\underset{i=1}{\cup}}\{p_{_{i}}\}$, the following Ostrowski-type inequality holds:

{\footnotesize
\Eq{*}{
\left|f(p)-\dfrac{1}{b-a}\int_{a}^{b}f(t)\,dt\right|
\leq \dfrac{1}{2}
&\max
\Bigg\{\Big\| f '\Big\|_a^{{p_{_{1}}}}({p_{_{1}}}-a)\,,\cdots\,, \Big\| f '\Big\|_{p_{_{i-1}}}^{p}(p-p_{_{i-1}})\,,\,\Big\| f '\Big\|_{p}^{p_{_{i}}}(p_{_{i}}-p)\,, \\
&\cdots\,,\, \Big\| f '\Big\|_{p_{_{n}}}^{b}(b-{p_{_{n}}})
\Bigg\}
+\max \Big\{f(a)+\overset{n}{\underset{i=1}{\sum}}{f(p_{_{i}})}\,,\,-\overset{n}{\underset{i=1}{\sum}}{f(p_{_{i}})}-f(b)\Big\}\,.
}
}
\\
Also, we investigate the possibility of proposing a refinement for Ostrowski inequality. We prove that if $f'\in L^{\infty}[a, b]$, then for any $p\in\,]a,b[$, we can restructure the inequality as follows:

\Eq{*}{
\Bigg|f(p)-\dfrac{1}{b-a}\int_{a}^{b}f(t)\,dt\Bigg|
\leq\min &\Bigg\{\left[\dfrac{1}{4}+\Bigg(\dfrac{p-\frac{a+b}{2}}{b-a}\Bigg)^2\right](b-a)\Big\|f' \Big \|_a^b \,,\,\\
&\,\,+\dfrac{1}{2}\,\max\bigg\{ (p-a)\Big\| f '\Big\|_{a}^{p}\,,\,(b-p)\Big\| f '\Big\|_{p}^{b}\bigg\}\Bigg\}.
}
The motivation, research background, methodology and crucial details are discussed in the introduction section.

\end{abstract}

\maketitle

\section*{Introduction} 
Throughout this paper $\N$, $\R$, and $\R_+$ denote the sets of natural, real, and non-negative numbers respectively. Unless or otherwise specified, the symbol $\Big\| f' \Big\|_a^b$ will be used to refer the supremum of the absolute value for the derivative of the differentiable function $f:[a,b]\to\R.$\\

In the year 1938, Swiss mathematician Alexander Ostrowski discovered that
if a function is differentiable and its derivative is bounded, then the difference between the function value at a point and the average value of the function over an interval can be bounded by a term involving the length of the interval and the bound on the derivative. Mathematically this can be represented as follows (see \cite{Ostrowski}): \\

"If $f:[a,b]\to\R$ is differentiable and $f'\in L^{\infty}[a,b]$, then for any $p\in\,]a,b[\,$, the following inequality holds:
\Eq{50}{
\Bigg|f(p)-\dfrac{1}{b-a}\int_{a}^{b}f(t)\,dt\Bigg|\leq\left[\dfrac{1}{4}+\Bigg(\dfrac{p-\frac{a+b}{2}}{b-a}\Bigg)^2\right](b-a)\Big\| f' \Big\|_a^b\,.\,^{^{^{^{^{"}}}}}
}
It is not difficult to validate that the rightmost expression of \eq{50} and the rightmost part of Ostrowski's inequality introduced the abstract are actually equivalent. Mathematicians have extensively investigated Ostrowski's inequality on various classes of functions. Some of the results are as follows:\\

"If $f\in C^n[a,b]$ and $f^{(n)}\in L^{\infty}[a, b]$, then for any $p\in\,]a,b[\,$, the following functional inequality is satisfied:
\Eq{*}{
\left| 
f(p) 
- \dfrac{1}{b - a} \int_a^b f(t) \, dt 
-\dfrac{1}{b - a}\sum_{k=1}^{n} \frac{(p-\frac{a+b}{2})^k}{k!} \int_a^b f^{(k)}(t) \, dt
\right|
\le 
\frac{(b - a)^n}{2^{n+1} n!} \, \Big\| f^{(n)} \Big\|_{a}^{b}\,.\,\, ^{^{^{^{^{^{"}}}}}}
}
A detailed proof of this theorem and other important results can be found in the paper \cite{derivative}.
\\

"If $f:[a,b]\to\R$ is convex, then for any $p\in\,]a,b[\,$, the following inequalities are satisfied:
\Eq{*}{
\dfrac{1}{b-a}\int_{a}^{b}f(t)\,dt-f(p)&\leq\dfrac{(p-a)(b-p)}{2(b-a)}\Bigg(\dfrac{f(b)-f(a)}{b-a}\Bigg)\\
&\mbox{and}\\
\qquad f(p)-\dfrac{1}{b-a}\int_{a}^{b}f(t)\,dt\leq\dfrac{1}{2(b-a)}\Bigg[(b-p)^2&\dfrac{f(p)-f(a)}{p-a}+(p-a)^2\dfrac{f(b)-f(p)}{b-p}\Bigg]\,.\, ^{^{^{^{^{^{"}}}}}}
}
These two inequalities together give a sandwich estimate of the deviation between  $f(p)$ and the integral mean of $f$ over $[a,b]$. The studies in \cite{Convex,sconvex} offer insights into this topic.\\

"If $f:[a,b]\to\R$ is a function of bounded variation, then for any $p\in\,]a,b[\,$, the following functional inequality holds:
\Eq{*}{
\Bigg|f(p)-\dfrac{1}{b-a}\int_{a}^{b}f(t)\,dt\Bigg|\leq 
\left[\dfrac{1}{2}+\left|\dfrac{p-\frac{a+b}{2}}{b-a}\right|\right]\bigvee_{a}^{b}(f).
^{^{^{^{^{^{"}}}}}}
}
where $\overset{b}{\underset{a}{\bigvee}}(f)$ denotes the total variation of $f$ in $[a,b]$. Regarding this more details can be found in the paper \cite{Bounded,Boundedd}.\\

"If $f:[a,b]\to\R$ is a absolutely continuous function, then for any $p\in\,]a,b[\,$, the following inequality is satisfied:

\Eq{*}{
\left| f(p) - \frac{1}{b - a} \int_a^b f(t) \, dt \right| \leq \left[ \frac{(p-a)(b - p)}{b - a} \right]{\Big\|f '\Big\|_a^b}\,;\quad
\mbox{where}\quad \Big\| f '\Big\|_a^b = \mbox{ess}\underset{t \in [a,b]}{\sup} |f'(t)|. 
^{^{^{^{"}}}}}
Relevant findings can be found in the works of \cite{Continuous, Absolutely}.\\

Besides these, under some strict conditions, several generalized and weighted versions have been formulated. Researchers also investigated relationship of this inequality with the theory of discrete means and numerical integration. Some of the works in this direction can be found in the papers \cite{Fink, George, Book, Barnett, Pecaric, Means, Ravi, Tail, pachpattee} and the references therein.\\

The paper can be divided into two parts. We start our investigation by deriving an Ostrowski-type inequality for the class of continuous functions that are not differentiable at a finite number of points. Besides the standard form, the assumption of differentiability is common for many other versions of Ostrowski-type inequalities. Although this regularity property is relaxed for the absolutely continuous functions and the functions with bounded variation, both these classes carry their own drawbacks. The class of absolutely continuous functions is a subclass of continuous functions. Similarly, the obtained upper bounds for the functions of bounded variations can be immensely large. Not to mention computing the total variation can be a difficult task too. We show that if $f:[a,b]\to\R$ is differentiable every where except at the points 
${p_{_{1}}},\,\cdots,\,{p_{_{n}}}$, then we can derive the following functional inequality:

\Eq{*}{
\left|f(p)-\dfrac{1}{b-a}\int_{a}^{b}f(t)\,dt\right|
\leq \dfrac{1}{2}
&\max
\Bigg\{\Big\| f '\Big\|_a^{{p_{_{1}}}}({p_{_{1}}}-a)\,,\cdots\,,\, \Big\| f '\Big\|_{p_{_{i-1}}}^{p}(p-p_{_{i-1}})\,,\Big\| f '\Big\|_{p}^{p_{_{i}}}(p_{_{i}}-p)\,, \\
&\qquad \cdots\cdots\cdots\,,\,\Big\| f '\Big\|_{p_{_{n}}}^{b}(b-{p_{_{n}}})
\Bigg\}+\max\Big\{f(a)-S\,,\,S-f(b)\Big\} \\
&\quad \Bigg(  \mbox{where}\quad S=f(a)+\overset{n}{\underset{i=1}{\sum}}{f(p_{_{i}})}
+f(b)\Bigg).
}
To reflect the importance of this result, some corollaries are also presented.\\

Ostrowski inequality is sharp which makes it challenging to come up with refinement options. However, some related works can be seen in the papers \cite{Ahmad, Liu, Refinement}. Under the same assumptions as those of the standard Ostrowski's inequality, a refinement of \eq{50} can be established as follows:
{\footnotesize
\Eq{*}{
\Bigg|f(p)-\dfrac{1}{b-a}\int_{a}^{b}f(t)\,dt\Bigg|\leq\min\left\{\left[\dfrac{1}{4}+\Bigg(\dfrac{p-\frac{a+b}{2}}{b-a}\Bigg)^2\right](b-a)\Big\|f \Big \|_a^b\,,\,\dfrac{1}{2}\,\max\bigg\{ (p-a)\Big\| f '\Big\|_{a}^{p}\,,\,(b-p)\Big\| f '\Big\|_{p}^{b}\bigg\}\right\}.}
}
Along with this result, a corollary and some inequalities associated with classical means are also included to demonstrate the impact of this finding.\\

We start our study with a discrete inequality.

\section{Main Results}
The following discrete inequality plays an important role throughout this paper. The result also reflects how effectively the concepts of discrete means can be applied in the theory of functional mean. The lemma can also be established by using mathematical induction. However here we give a short and simple proof.

\Lem{22}{Let $a_1,\cdots,a_n\in\R$ and $b_1,\cdots,b_n\in\R_+$, then the following discrete functional inequality is satisfied
\Eq{1}{
\min\Bigg(\dfrac{a_1}{b_1},\cdots,\dfrac{a_n}{b_n}
\Bigg)\leq 
\dfrac{a_1+\cdots+a_n}{b_1+\cdots+b_n}
\leq\max\Bigg(\dfrac{a_1}{b_1},\cdots,\dfrac{a_n}{b_n}
\Bigg).
}
}
\begin{proof}
The expression $\displaystyle{\dfrac{a_1+\cdots+a_n}{b_1+\cdots+b_n}}$ can be re-written as the following convex combination
\Eq{*}{
\dfrac{b_1}{b_1+\cdots+b_n}\Bigg(\dfrac{a_1}{b_1}\Bigg)+\cdots+\dfrac{b_n}{b_1+\cdots+b_n}\Bigg(\dfrac{a_n}{b_n}\Bigg).
}
Hence by the mean property, inequality \eq{1} is obvious.
\end{proof}
In the next result, we develop an Ostrowski-type inequality for a broader class of functions.
\Thm{88}{Let $f\in C[a,b]$ is differentiable everywhere except at the points ${p_{_{1}}},\,\cdots\,,{p_{_{n}}}$. Then for any $p\in\,]a,b[\,\setminus\,\{p_1,\cdots,p_n\}$, the following Ostrowski-type inequality is satisfied:
\Eq{907}{
\left|f(p)-\dfrac{1}{b-a}\int_{a}^{b}f(t)\,dt\right|
\leq \dfrac{1}{2}
&\max
\Bigg\{\Big\| f '\Big\|_a^{{p_{_{1}}}}({p_{_{1}}}-a)\,,\cdots\,,\, \Big\| f '\Big\|_{p_{_{i-1}}}^{p}(p-p_{_{i-1}})\,,\Big\| f '\Big\|_{p}^{p_{_{i}}}(p_{_{i}}-p)\,, \\
&\qquad \cdots\cdots\cdots\,,\,\Big\| f '\Big\|_{p_{_{n}}}^{b}(b-{p_{_{n}}})
\Bigg\}+\max\Big\{f(a)-S\,,\,S-f(b)\Big\} \\
&\quad \Bigg(  \mbox{where}\quad S=f(a)+\overset{n}{\underset{i=1}{\sum}}{f(p_{_{i}})}
+f(b)\Bigg).
}
}
\begin{proof}
To prove theorem, We proceed as follows:
\Eq{*}{
f(p)-\dfrac{1}{b-a}\int_{a}^{b}f(t)\,dt &=\dfrac{1}{b-a}\int_{a}^{b}\Big(f(p)-f(t)\Big)\,dt\\
&= \dfrac{1}{b-a}\int_{a}^{b}\Big(f(p)-f(t)+S-f(a)\Big)\,dt+f(a)-S\,.
}
We consider the partition of the interval $[a,b]$ as follows:
$$\P:=\Big\{a=p_0,\cdots,p_{_{i-1}},p,p_{_{i}},\cdots, p_{_{n+1}}=b\Big\}\,.$$
Now to extend the rightmost part of the above equality. First we consider partition $\P$ and use the additivity property of integration at the partitioning points. After that, we utilize \eq{1} and subinterval-wise implement the Cauchy mean value theorem to apply the $L^{\infty}$ norms. The computation is as follows:
\Eq{*}{
&\dfrac{\displaystyle \begin{array}{c}
\displaystyle \int_{p_{_{0}}}^{p_{_{1}}}\big(f(p_{_{1}})-f(t)\big)\,dt+\cdots+\int_{p_{_{i-1}}}^{p}\big(f(p)-f(t)\big)\,dt \\
\displaystyle\qquad\qquad 
+\int_{p}^{p_{_{i}}}\big(f({p_{_{i}}})-f(t)\big)\,dt+\cdots+
\int_{p_{_{n}}}^{p_{_{n+1}}}\big(f(p_{_{n+1}})-f(t)\big)\,dt
\end{array}}
{\displaystyle 
 \big(p_{_{1}}-p_{_{0}}\big)+\cdots +\big(p-p_{_{i-1}}\big)+\big(p_{_{i}}-p\big)+\cdots+\big(p_{_{n+1}}-p_{_n} \big)}+f(a)-S\\
&\leq \max\Bigg\{\dfrac{1}{p_{_{1}}-p_{_{0}}}\int_{p_{_{0}}}^{p_{_{1}}}\big(f(p_{_{1}})-f(t)\big)\,dt,\cdots,\dfrac{1}{p-p_{_{i-1}}}\int_{p_{_{i-1}}}^{p}\big(f(p)-f({z})\big)\,dt\\
&\qquad\quad\quad \dfrac{1}{p_{_{i}}-p}\int_{p}^{p_{_{i+1}}}\big(f(p)-f({t})\big)\,dt,\cdots, \dfrac{1}{p_{_{n+1}}-p_{_{n}}}\int_{p_{_{n}}}^{p_{_{n+1}}}\big(f({p_{_{n+1}}})-f({t})\big)\,dt\Bigg\}
+f(a)-S\\
&\leq\max\Bigg\{\dfrac{\Big\| f '\Big\|_{p_{_{0}}}^{p_{_{1}}}}{p_{_{1}}-p_{_{0}}}\int_{p_{_{0}}}^{p_{_{1}}}(p_{_{1}}-t)\,dt,\cdots,\dfrac{\Big\| f '\Big\|_{p_{_{i-1}}}^{p}}{p-p_{_{i-1}}}\int_{p_{_{i-1}}}^{p}(p-t)\,dt,\\
&\qquad\quad\quad \dfrac{\Big\| f '\Big\|_{p}^{p_{_{i}}}}{p_{_{i}}-p}\int_{p_{_{i}}}^{p}(p_{_{i}}-t)\,dt,\cdots,
\dfrac{\Big\| f '\Big\|_{p_{_{n}}}^{p_{_{n+1}}}}{p_{_{n+1}}-p_{_{n}}}\int_{p_{_{n}}}^{p_{_{n+1}}}(p_{_{n+1}}-t)\,dt
\Bigg\}+f(a)-S\\
&=\dfrac{1}{2}\max\Bigg\{\Big\| f '\Big\|_{p_{_{0}}}^{p_{_{1}}}(p_{_{1}}-p_{_{0}})\,, \cdots\cdots \cdots,\Big\| f '\Big\|_{p_{_{i-1}}}^{p}(p-p_{_{i-1}})\,,\\
&\qquad \qquad\qquad\qquad 
\Big\| f '\Big\|_{_{p}}^{p_{_{i}}}(p_{_{i}}-p)\,,\cdots\cdots\cdots,\Big\| f '\Big\|_{p_{_{n}}}^{p_{_{n+1}}}(p_{_{n+1}}-p_{_{n}})\Bigg \}
+f(a)-S\,.
}
Therefore, we obtain the following inequality
\Eq{900}{ f(p)-\dfrac{1}{b-a}\int_{a}^{b}f(t)\,dt 
&\leq \dfrac{1}{2}\max\Bigg\{\Big\| f '\Big\|_{a}^{p_{_{1}}}(p_{_{1}}-a)\,, \cdots\cdots ,\Big\| f '\Big\|_{p_{_{i-1}}}^{p}(p-p_{_{i-1}})\,,\\
&\qquad \qquad \,\,\,\,\,
\Big\| f '\Big\|_{_{p}}^{p_{_{i}}}(p_{_{i}}-p)\,,\cdots\cdots,\Big\| f '\Big\|_{p_{_{n}}}^{b}(b-p_{_{n}})\Bigg \}
+f(a)-S\,.
}
We can also compute the following inequality
\Eq{*}{
\dfrac{1}{b-a}\int_{a}^{b}f(t)\,dt-f(p) &=\dfrac{1}{b-a}\int_{a}^{b}\Big(f(t)-f(p)\Big)\,dt\\
&= \dfrac{1}{b-a}\int_{a}^{b}\Big(f(t)-f(p)+f(b)-S\Big)\,dt+S-f(b).}
Proceeding with similar mathematical techniques as before, we can extend the rightmost part of above integral equality as follows
\Eq{*}{
&\dfrac{\displaystyle \begin{array}{c}
\displaystyle \int_{p_{_{0}}}^{p_{_{1}}}\big(f(t)-f(p_{_{0}})\big)\,dt+\cdots+\int_{p_{_{i-1}}}^{p}\big(f(t)-f(p_{_{i-1}})\big)\,dt\\
\displaystyle \qquad\qquad+\int_{p}^{p_{_{i}}}\big(f(t)-f(p)\big)\,dt+\cdots+
\int_{p_{_{n}}}^{p_{_{n+1}}}\big(f(t)-f({p_{_{n}}})\big)\,dt
\end{array}}
{\displaystyle{\big(p_{_{1}}-p_{_{0}}\big)+\cdots \big(p-p_{_{i-1}}\big)+\big(p_{_{i}}-p_{_{i-1}}\big)+\cdots\big(p_{_{n+1}}-p_{_{n}} \big)} }+S-f(b)\\
&\leq \max\Bigg\{\dfrac{1}{p_{_{1}}-p_{_{0}}}\int_{p_{_{0}}}^{p_{_{1}}}\big(f(t)-f(p_{_{0}})\big)\,dt,\cdots, \dfrac{1}{p-p_{_{i-1}}}\int_{p_{_{i-1}}}^{p}\big(f(t)-f(p_{_{i-1}})\big)\,dt\\
&\qquad\quad\quad \dfrac{1}{p_{_{i}}-p}\int_{p}^{p_{_{i+1}}}\big(f({t})-f(p)\big)\,dt,\cdots, \dfrac{1}{p_{_{n}}-p_{_{n+1}}}\int_{p_{_{n}}}^{p_{_{n+1}}}\big(f({t})-f(p_{_{n+1}})\big)\,dt\Bigg\}+S-f(b)\\
&\leq\max\Bigg\{\dfrac{\Big\| f '\Big\|_{p_{_{0}}}^{p_{_{1}}}}{p_{_{1}}-p_{_{0}}}\int_{p_{_{0}}}^{p_{_{1}}}(t-p_{_0})\,dt,\cdots,\dfrac{\Big\| f '\Big\|_{p_{_{i-1}}}^{p}}{p-p_{_{i-1}}}\int_{p_{_{i-1}}}^{p}(t-p_{_{i-1}})\,dt,\\
&\qquad\quad\quad \dfrac{\Big\| f '\Big\|_{p}^{p_{_{i}}}}{p_{_{i}}-p}\int_{p_{_{i}}}^{p}(t-p)\,dt,\cdots,
\dfrac{\Big\| f '\Big\|_{p_{_{n}}}^{p_{_{n+1}}}}{p_{_{n+1}}-p_{_{n}}}\int_{p_{_{n}}}^{p_{_{n+1}}}(t-p_{_{n}})\,dt
\Bigg\}+S-f(b)\\
&=\dfrac{1}{2}\max\Bigg\{\Big\| f '\Big\|_{p_{_0}}^{p_{_{1}}}(p_{_{1}}-p_{_{0}})\,, \cdots\cdots \cdots,\Big\| f '\Big\|_{p_{_{i-1}}}^{p}(p-p_{_{i-1}})\,,\\
&\qquad \qquad\qquad\qquad 
\Big\| f '\Big\|_{p}^{p_{_{i}}}(p_{_{i}}-p)\,,\cdots\cdots\cdots,\Big\| f '\Big\|_{p_{_{n}}}^{p_{_{n+1}}}(p_{_{n+1}}-p_{_{n}})\Bigg \}+S-f(b)\,.
}
Thus we have
\Eq{901}{\dfrac{1}{b-a}\int_{a}^{b}f(t)\,dt -f(p)
&\leq \dfrac{1}{2}\max\Bigg\{\Big\| f '\Big\|_a^{p_{_{1}}}(p_{_{1}}-a)\,, \cdots\cdots ,\Big\| f '\Big\|_{p_{_{i-1}}}^{p}(p-p_{_{i-1}})\,,\\
&\qquad \qquad \,\,\,\,\,
\Big\| f '\Big\|_{_{p}}^{p_{_{i}}}(p_{_{i}}-p)\,,\cdots\cdots,\Big\| f '\Big\|_{p_{_{n}}}^b(b-p_{_{n}})\Bigg \}
+S-f(b)\,.
}
Combining the inequalities \eq{900} and \eq{901} we achieve \eq{907}. This completes the proof.
\end{proof}
In the next result, instead of a differentiable point, we consider a point of non-differentiability and obtain an Ostrowski-type inequality.
\Thm{89}{Let $f\in C[a,b]$ is non-differentiable only at the points ${p_{_{1}}},\cdots, \, {p_{_{n}}}\in\,]a,b[$. Then for any ${p_{_{i}}}$ $(i\in\{1,\cdots,n\})$ the following Ostrowski's type inequality is satisfied.
\Eq{*}{
\left|f({p_{_{i}}})-\dfrac{1}{b-a}\int_{a}^{b}f(t)\,dt\right|\leq \dfrac{1}{2}&\max
\Bigg\{\Big\| f '\Big\|_{p_{_{i-1}}}^{p_{_{i}}}(p_{_{i}}-p_{_{i-1}})\,,\,i\in\{1,\cdots ,n\}
\Bigg\}\\
+&\max\Big\{f(a)+f(p_{_{i}})-S\,,\,S-f(p_{_{i}})-f(b)\Big\}\,\\
&\qquad \Bigg(  \mbox{where}\quad S=f(a)+\overset{n}{\underset{i=1}{\sum}}{f(p_{_{i}})}+f(b)
\Bigg)\,.
}
}
\begin{proof}
The proof of this theorem is similar to \thm{88}. Hence we decided to leave it to the readers.
\end{proof}
Both of the above results can also be generalized to countably infinite number of discontinuity provided the term $S={\underset{i\in\N}{\sum}}{f(p_{_{i}})}$ is finite. Next we discuss the concept of refinement in inequalities.\\

For the classical Ostrowski's inequality, we can propose several versions of refinement. For instance, we can figure out a function $\varphi:[a,b]\to\R$ such that it will satisfy the following inequality under the required pre-assumed conditions on $f$
\Eq{*}{
\Bigg|f(p)-\dfrac{1}{b-a}\int_{a}^{b}f(t)\,dt\Bigg|
\leq \varphi(p)
\leq\left[\dfrac{1}{4}+\Bigg(\dfrac{p-\frac{a+b}{2}}{b-a}\Bigg)^2\right](b-a)\Big\| f' \Big\|_a^b\,.
}
Another way of improving Ostrowski's inequality is to figure out a meaningful function $\psi:[a,b]\to\R$ and represent the inequality as below
\Eq{*}{
\Bigg|f(p)-\dfrac{1}{b-a}\int_{a}^{b}f(t)\,dt\Bigg|
\leq\min\Bigg\{\psi(p)\,,\left[\dfrac{1}{4}+\Bigg(\dfrac{p-\frac{a+b}{2}}{b-a}\Bigg)^2\right](b-a)\Big\| f' \Big\|_a^b\Bigg\}\,.
}
In the next result, we propose a refinement of Ostrowski's inequality.
\Thm{2}{
Let $f\in C[a,b]$ and $f'\in L^{\infty}[a,b].$ Then for any $p\in\,]a,b[\,$, the following inequality holds:
{\footnotesize
\Eq{100}{
\Bigg|f(p)-\dfrac{1}{b-a}\int_{a}^{b}f(t)\,dt\Bigg|
\leq\min\left\{\left[\dfrac{1}{4}+\Bigg(\dfrac{p-\frac{a+b}{2}}{b-a}\Bigg)^2\right](b-a)\Big\| f' \Big\|_a^b\,,\,\dfrac{1}{2}\,\max\bigg\{ (p-a)\Big\| f '\Big\|_{p_{_{n}}}^{p_{_{n+1}}}\,,\,(b-p)\Big\| f '\Big\|_{p}^{b}\bigg\}\right\}\,.
}}
}
\begin{proof}
The expression $\displaystyle{\dfrac{1}{b-a}\int_{a}^{b}f(t)\,dt}$ can be re-written as follows:
$$\dfrac{1}{b-a}\int_{a}^{b}f(t)\,dt:=\frac{\int_{a}^{p}f(t)\,dt+\int_{p}^{b}f(t)\,dt}{(p-a)+(b-p)}.$$
Therefore, by applying the inequality \eq{1} for $n=2$ here, we can conclude the following
\begin{small}
\Eq{*}{
\min\bigg\{\dfrac{1}{p-a}\int_{a}^{p}f(t)\,dt\,,\, \dfrac{1}{b-p}\int_{p}^{b}f(t)\,dt\bigg\}
\leq\dfrac{1}{b-a}\int_{a}^{b}f(t)\,dt\leq\max\bigg\{\dfrac{1}{p-a}\int_{a}^{p}f(t)\,dt\,,\, \dfrac{1}{b-p}\int_{p}^{b}f(t)\,dt\bigg\}.
}
\end{small}
Without loss of generality, we assume that
\Eq{*}{
\dfrac{1}{p-a}\int_{a}^{p}f(t)\,dt\leq\dfrac{1}{b-a}\int_{a}^{b}f(t)\,dt\leq \dfrac{1}{b-p}\int_{p}^{b}f(t)\,dt\,.
}
Proceeding with this assumption and later using the mean value theorem as well as the $L^{\infty}$ norm, we obtain the following inequality
\Eq{77}{
\dfrac{1}{b-a}\int_{a}^{b}f(t)\,dt-f(p)&\leq \dfrac{1}{b-p}\int_{p}^{b}f(t)\,dt-f(p)\\
&=\dfrac{1}{b-p}\int_{p}^{b}\Big(f(t)-f(p)\Big)\,dt
=\dfrac{1}{b-p}\int_{p}^{b}\dfrac{f(t)-f(p)}{t-p}(t-p)\,dt\\
&\qquad\qquad\qquad\qquad\qquad\qquad \quad
\leq \dfrac{\Big\| f '\Big\|_{p}^{b}}{b-p}\int_{p}^{b}(t-p)\,dt=\dfrac{1}{2}\,
{\Big\| f '\Big\|_{p}^{b}}(b-p)\,.
}
Similarly, under our assumptions, we can show that 
\Eq{78}{
f(p)-\dfrac{1}{b-a}\int_{a}^{b}f(t)\,dt\leq f(p)-\dfrac{1}{p-a}\int_{a}^{p}f(t)\,dt \leq \dfrac{1}{2}\,{\Big\| f '\Big\|_{a}^{p}}(p-a)\,.
}
Combining the two inequalities \eq{77} and \eq{78}, we can conclude the following
\Eq{*}{
\Bigg|f(p)-\dfrac{1}{b-a}\int_{a}^{b}f(t)\,dt\Bigg|\leq\dfrac{1}{2}\,\max\bigg\{ (p-a)\Big\| f '\Big\|_{a}^{p}\,,\,(b-p)\Big\| f '\Big\|_{p}^{b}\bigg\}\,.
}

This together with \eq{50} validates \eq{100} and completes the proof.
\end{proof}
Instead of considering only one interior point $p$ in \thm{2}, we can try to consider multiple interior points and follow the same methodology to obtain a another Ostrowski-type inequality. This ensures a better refinement.\\

The following result can be easily derived by using the same procedure as in \thm{2}. Hence, we decided to include the statement only. 

\Thm{22}{
If the function $f:[a,b]\to\R$ is continuous and differentiable except the point $p\in\,]a,b[\,$, then the following functional inequality holds
{\footnotesize
\Eq{2233}{
\Bigg|f(p)-\dfrac{1}{b-a}\int_{a}^{b}f(t)\,dt\Bigg|
\leq\min\left\{\left[\dfrac{1}{4}+\Bigg(\dfrac{p-\frac{a+b}{2}}{b-a}\Bigg)^2\right](b-a)\Big\|f \Big \|_a^b\,,\dfrac{1}{2}\,\max\bigg\{ (p-a)\Big\| f '\Big\|_{a}^{p}\,,\,(b-p)\Big\| f '\Big\|_{p}^{b}\bigg\}\right\}
\,.}
}
}

Now in the last part of our paper, we show the application of \thm{2}, we demonstrate several enhanced inequalities between various classical discrete means. First, we recall various notions and terminologies of these standard means
Let $a,b\in\R_+$. Then we define some of the important discrete means as follows:
\Eq{*}{
\mbox{ Arithmetic mean, } \A:&=\dfrac{a+b}{2}\\
\mbox{Geometric mean, } \G:&=\sqrt{ab}\\
\mbox{Harmonic mean, } \H:&=\dfrac{2ab}{a+b}\\
\mbox{Logarithmic mean, } \L:&=\dfrac{b-a}{ln(b)-ln(a)}
}
In a non-empty and non-singleton interval $[a,b]\in\R_+$, the order of these means are as follows
\Eq{556}{
a<\H<\G<\L<\A<b\,.
}
The following theorem contains the improved versions of several inequalities obtained by Dragomir and Wang in their papers \cite{meanss, Wang} showcasing relationships among various discrete means. Following the same procedure, One can work on the refinements of other inequalities mentioned in these papers that involves other special means. 
\Thm{7862}{
Let $a,b\in\R_+$ with $a<b$. Then the following inequalities are satisfied:
\Eq{557}{
(i)\quad 0&\leq \A-\L\leq \dfrac{\A\L(b-a)}{4a^2}\,,\\
(ii)\quad 0&\leq \L-\G\leq \min\Bigg\{ \dfrac{\G\L(b-a)}{a^2}\bigg[\dfrac{1}{4}+\bigg(\dfrac{\G-\A}{b-a}\bigg)^2\bigg]\,,\,\dfrac{\G-a}{2a^2}\Bigg\}\,,\\
(iii)\quad 0&\leq \L-\H\leq \min\Bigg\{ \dfrac{\H\L(b-a)}{a^2}\bigg[\dfrac{1}{4}+\bigg(\dfrac{\H-\A}{b-a}\bigg)^2\bigg]\,,\,\dfrac{b-\H}{2\H^2}\Bigg\}
.}
}
\begin{proof}
To prove these inequalities, we substitute $f(t):=1/t$ in the inequality \eq{2233} and obtain the following
\Eq{*}{
\left|\dfrac{\L-p}{\L p}\right|
\leq \min\Bigg\{\left(\dfrac{(p-a)^2+(b-p)^2}{2\cdot(b-a)}\right)\dfrac{1}{a^2}\,,\,\dfrac{1}{2}\,\max\bigg\{\dfrac{p-a}{a^2}\,,\,\dfrac{b-p}{p^2}\bigg\}\Bigg\}\,.
}
Now replacing $p$ with $\A$, $\G$, and $\H$ respectively and using the above inequality we can conclude \eq{557}. This completes the proof.
\end{proof}
There are several directions in which the findings of this paper can be further generalized. One natural extension is to relax the current assumption that discontinuities are allowed only at a finite number of points. By incorporating the notion of one-sided continuity and applying analogous analytical techniques, it may be possible to extend the results to broader classes of functions, thereby enhancing their applicability in more general settings.


\begin{thebibliography}{10}
\bibitem{Ostrowski}
Ostrowski, Alexander. 
\newblock{Über die Absolutabweichung einer differentiierbaren Funktion von ihrem Integralmittelwert.}
\newblock{Commentarii Mathematici Helvetici 10, no. 1 (1937): 226-227.}

\bibitem{derivative}
Cerone, Pietro, Sever S. Dragomir, and John Roumeliotis. 
\newblock{Some ostrowski type inequalities for $n$-time differentiable mappings and applications.}
\newblock{Demonstratio Mathematica 32, no. 4 (1999): 697-712.}


\bibitem{Convex}
Dragomir, Sever S. 
\newblock{Ostrowski’s Type Inequalities for Convex Functions and Applications.}
\newblock{RGMIA Research Report Collection 2, no. 1 (1999): Article 4.}
  
\bibitem{sconvex}
Alomari, Mohammad, Maslina Darus, Sever Silvestru Dragomir, and Pietro Cerone. \newblock{Ostrowski type inequalities for functions whose derivatives are s-convex in the second sense.}
\newblock{Applied mathematics letters 23, no. 9 (2010): 1071-1076.}
  
\bibitem{Bounded}  
Dragomir, S. S. 
\newblock{The Ostrowski integral inequality for mappings of bounded variation.}
\newblock{Bulletin of the Australian Mathematical Society 60, no. 3 (1999): 495-508.}

\bibitem{Boundedd}
Budak, Huseyin, and Mehmet Zeki Sarikaya. 
\newblock{On generalization of weighted Ostrowski type inequalities for functions of bounded variation.}
\newblock{Asian-European Journal of Mathematics (2015).}

\bibitem{Continuous}
Dragomir, Sever S. 
\newblock{Some companions of Ostrowski's inequality for absolutely continuous functions and applications.}
\newblock{Bulletin of the korean mathematical society 42, no. 2 (2005): 213-230.}

\bibitem{Absolutely}
Kikianty, Eder, Sever S. Dragomir, and Pietro Cerone. 
\newblock{Ostrowski type inequality for absolutely continuous functions on segments in linear spaces.}
\newblock{Research report collection 10, no. 3 (2007).}

\bibitem{Fink}
Fink, A. M. 
\newblock{Bounds on the deviation of a function from its averages.}
\newblock{Czechoslovak Mathematical Journal 42, no. 2 (1992): 289-310.}
\bibitem{George}
Anastassiou, George A. 
\newblock{Ostrowski type inequalities.}
\newblock{Proceedings of the American Mathematical Society 123.12 (1995): 3775-3781.}

\bibitem{Book}
Dragomir, Sever Silvestru, and Themistocles M. Rassias, eds. 
\newblock{Ostrowski type inequalities and applications in numerical integration.} \newblock{Dordrecht: Kluwer Academic, 2002.}

\bibitem{Barnett}
Barnett, Neil S., Sever Silvestru Dragomir, and Ian Gomm. 
\newblock{A companion for the Ostrowski and the generalised trapezoid inequalities.} \newblock{Mathematical and Computer Modelling 50, no. 1-2 (2009): 179-187.}

\bibitem{Pecaric}
Pečarić, Josip E. 
\newblock{Some further remarks on the Ostrowski generalization of Čebyšev's inequality.}
\newblock{Journal of mathematical analysis and applications 123, no. 1 (1987): 18-33.}

\bibitem{Means}
Tunç, Mevlüt. 
\newblock{Ostrowski-type inequalities via h-convex functions with applications to special means.}
\newblock{Journal of Inequalities and Applications 2013 (2013): 1-10.}

\bibitem{Ravi}
Agarwal, Ravi P., M-J. Luo, and R. K. Raina. 
\newblock{On Ostrowski type inequalities.}
\newblock{Fasciculi Mathematici (2016).}

\bibitem{Tail}
Acu, Ana Maria, Heiner Gonska, and I. Raşa. 
\newblock{Grüss-type and Ostrowski-type inequalities in approximation theory.}
\newblock{Ukrainian Mathematical Journal 63, no. 6 (2011): 843-864.}

\bibitem{pachpattee}
Pachpatte, B. G. 
\newblock{On an inequality of Ostrowski type in three independent variables.}
\newblock{Journal of Mathematical Analysis and Applications 249, no. 2 (2000): 583-591.}

\bibitem{Ahmad}
Ahmad, Hijaz, Muhammad Tariq, Soubhagya Kumar Sahoo, Sameh Askar, Ahmed E. Abouelregal, and Khaled Mohamed Khedher. 
\newblock{Refinements of Ostrowski type integral inequalities involving Atangana–Baleanu fractional integral operator.}
\newblock{Symmetry 13, no. 11 (2021): 2059.}


\bibitem{Liu}
Liu, Wenjun, and Yun Sun. 
\newblock{A Refinement of the Companion of Ostrowski inequality for functions of bounded variation and Applications.}
\newblock{arXiv preprint arXiv:1207.3861 (2012).}

\bibitem{Refinement}
Dragomir, Sever S. 
\newblock{Refinements of the Ostrowski inequality in terms of the cumulative variation and applications.}
\newblock{Analysis 34, no. 2 (2014): 223-240.}

\bibitem{meanss}
Dragomir, S. S., and S. Wang. 
\newblock{An inequality of Ostrowski-Gr\"uss' type and its applications to the estimation of error bounds for some special means and for some numerical quadrature rules.} 
\newblock{Computers and Mathematics with Applications 33, no. 11 (1997): 15-20.}

\bibitem{Wang}
Dragomir, S. S., and S. Wang. 
\newblock{Applications of Ostrowski's inequality to the estimation of error bounds for some special means and for some numerical quadrature rules.}
\newblock{Applied Mathematics Letters 11, no. 1 (1998): 105-109.}
\end{thebibliography}
\bibliographystyle{plain}

\end{document}